\documentclass{amsart}

\usepackage{graphics}
\usepackage{epsfig}
\usepackage{amsmath}
\usepackage{amssymb,fancyhdr,txfonts,pxfonts}

\usepackage{paralist}

\usepackage{pstricks,pst-text,pst-grad,pst-node,pst-3dplot,pstricks-add,pst-poly,pst-coil} 
\usepackage{pst-fun,pst-blur} 

\usepackage{hyperref}

\newtheorem{theorem}{Theorem}
\theoremstyle{plain}

\newtheorem{corollary}{Corollary}

\newtheorem{definition}{Definition}

\newtheorem{lemma}{Lemma}

\newtheorem{proposition}{Proposition}

\begin{document}

\title[Endomorphisms Algebra of Translations Group in Affine Plane]{The Endomorphisms Algebra of Translations Group
and Associative Unitary Ring of Trace-Preserving Endomorphisms in Affine Plane}

\author[Orgest ZAKA]{Orgest ZAKA}
\address{Orgest ZAKA: Department of Mathematics-Informatics, Faculty of Economy and Agribusiness, 
Agricultural University of Tirana, Tirana, Albania.}
\email{ozaka@ubt.edu.al, gertizaka@yahoo.com, ozaka@risat.org}

\dedicatory{Dedicated to David Hilbert and Emil Artin}


\begin{abstract}
This paper introduces a description of Endomorphisms of the translation group in an affine plane,  will define the addition and composition of the set of endomorphisms and specify the neutral elements associated with these two actions and present the Endomorphism algebra thereof  
will distinguish the Trace-preserving endomorphism algebra in affine plane, and prove that the set of Trace-preserving endomorphism associated with the 'addition' action forms a commutative group. We also try to prove that the set of trace-preserving endomorphism, together with the two actions, in it, 'addition' and 'composition' forms an associative and unitary ring.
\end{abstract}


\maketitle

\section{Introduction}

The foundations for the study of the connections between axiomatic geometry and algebraic structures were set forth by D. Hilbert \cite{Hilbert1959geometry}, recently elaborated and extended in terms of the algebra of affine planes in, for example, \cite{Kryftis2015thesis}, [3, §IX.3, p.574]. Also great contributions in this direction have been made by, E.Artin in \cite{11},  H. S. M. Coxeter, in \cite{14},  Marcel Berger in \cite{Marcel.Berger.Geometry1.and.Geometry2}, Robin Hartshrone in \cite{Hartshorne.Foundations}, etc. Even earlier, in my works '\cite{1}, \cite{2}, \cite{3}, \cite{4}, \cite{5}, \cite{6}, \cite{7}, \cite{15}, \cite{Zaka_2016}' I have brought up quite a few interesting facts about the association of algebraic structures with affine planes and with 'Desargues affine planes', and vice versa. 

In this paper, we will present a description of Endomorphisms of the translation group in affine plane. We will define the addition and composition of the set of endomorphisms of the translation group in the affine plane. We will specify the neutral elements associated with these two actions, which we will call, 'zero endomorphism' and 'unitary endomorphism', and present the Endomorphism algebra of translations in affine plane. We will distinguish, as a substructure of Endomorphisms algebra, the Trace-preserving endomorphism algebra in affine planes, and prove that the set of Trace-preserving endomorphism associated with the 'addition' action forms a commutative group, also prove that the set of trace-preserving endomorphism, together with the two actions, in it, 'addition' and 'composition' forms an associative and unitary ring.

In paper \cite{5}, we have done a detailed, careful description, of translations and dilation’s in affine planes. We have proven that the set of dilation’s regarding the composition action is a group, and set of translations is a commutative group. We have proved that translation group is a normal subgroup of the group of dilation’s. We have described and specified what we would call a translation or dilation trace, we have also defined the direction of an translation as an equivalence class of parallelism in affine planes. This will help us in this paper, as trace-preserving endomorphisms to retain those endomorphisms that operate on these equivalence classes according to parallelism.

\section{Preliminaries}
Let $\mathcal{P}$ be a nonempty set, which is called set of points, $\mathcal{L}$ a nonempty subset of $\mathcal{P}$, which is called set of lines, and an incidence relation $\mathcal{I}\subseteq \mathcal{P} \times \mathcal{L}$.  For a point $P \in \mathcal{P}$ and an line $\ell \in \mathcal{L}$, the fact $(P,\ell )\in \mathcal{I}$, (equivalent to $P\mathcal{I}\ell $) we mark $P\in \ell$ and read \textit{point} $P$ \textit{ is incident with a line} $\ell$ or a line passes through points $P$ (contains point $P$).

\begin{definition} \label{def.affine.plane} 
\cite{11},\cite{13},\cite{3} Affine plane is called the incidence
structure  $\mathcal{A}=(\mathcal{P},\mathcal{L},\mathcal{I})$ that satisfies
the following axioms:
\begin{description}
	\item[A.1] For every two different points $P$ and $Q$ $\in $ $\mathcal{P}$, there
exists exactly one line $\ell \in $ $\mathcal{L}$ incident with that points.
  \item[A.2] For a point $P$ $\in $ $\mathcal{P}$, and an line $\ell \in $ $\mathcal{L}$
such that $(P,\ell )\notin $ $\mathcal{I}$, there exists one and only one
line $r\in $ $\mathcal{L}$, incident with the point $P$ and such that $\ell
\cap r=\varnothing .$
  \item[A.2] In $\mathcal{A}$ there are three non-incident points with a line.
\end{description}
\end{definition}

Whereas a line of the affine plane we consider as sets of points of affine
plane with her incidents. Axiom A.1 implicates that tow different lines of $
\mathcal{L}$ many have a common point, in other words tow different lines of 
$\mathcal{L}$ either have no common point or have only one common point.

\begin{definition}
Two lines $ \ell , m \in \mathcal{L}$  that are matching or do not have
common point are called \textbf{parallel} and in this case is written  $\ell
\parallel m,$ and when they have only one common point we say that they are
expected.
\end{definition}

\begin{definition} 
Let it be $\mathcal{A}=(\mathcal{P},\mathcal{L},\mathcal{I})$ an
affine plane and $\mathcal{S}$=$\{\psi :\mathcal{P}\rightarrow \mathcal{P}%
| $ where $\psi -$is bijection$\}$  set of bijections to set points $\mathcal{P}$ on yourself. Collineation of affine plane $\mathcal{A}$, called a bijection $\psi \in \mathcal{S}$, such that%
\begin{equation}
\forall \ell \in \mathcal{L},\psi \left( \ell \right) \in \mathcal{L},
\label{eq}
\end{equation}
\end{definition}

Otherwise, a collineation of the affine plane $\mathcal{A}$ is a bijection
of set $\mathcal{P}$ on yourself, that preserves lines.  It is known
that the set of bijections to a set over itself is a group on associated
with the binary action '$\circ$' of composition in it, which is known as
total group or symmetric groups.

\begin{definition} \label{def.fixpoint} \cite{5}
An point $P$ of the affine plane $\mathcal{A}$ called fixed point his
associated with a collineation $\delta ,$ if coincides with the image itself
 $\delta (P),$ briefly when, $P=\delta (P).$
\end{definition}

\begin{definition} \label{def.dilation} \cite{11}, \cite{5}
A \textbf{Dilation} of an affine plane $\mathcal{A}=(\mathcal{P},%
\mathcal{L},\mathcal{I})$ called a its collineation $\delta $ such that%
\begin{equation}
\forall P\neq Q\in \mathcal{P},\delta \left( PQ\right) \Vert PQ
\end{equation}
\end{definition}

Let it be $\mathbf{Dil}_{\mathcal{A}}=\left\{ \delta \in \mathbf{Col}_{%
\mathcal{A}}|\delta -\text{is a dilation of }\mathcal{A}\right\} $ the
dilation set of affine plane $\mathcal{A}=(\mathcal{P},\mathcal{L},%
\mathcal{I})$.

\begin{theorem} \cite{5}, \cite{15}
The dilation set $\mathbf{Dil}_{\mathcal{A}}$ \ of affine plane $\mathcal{A%
}$ forms a group with respect to composition $\circ .$
\end{theorem}
\proof See, Theorem 2.4, in \cite{5}. \qed

\begin{definition} \label{def.trace} 
Let it be $\delta $ an dilation of affine plane $\mathcal{A}=(\mathcal{P},%
\mathcal{L},\mathcal{I})$, and $P$ a point in it. Lines that passes by points $P$ and $%
\delta (P),$ called \textbf{trace} of points $P$ regarding dilation $\delta .$
\end{definition}

Every point of a traces of a not-fixed point, to an affine plane
associated with its dilation has its own image associated with that dilation
in the same traces (see \cite{5}). We also know the result: If an affine plane $\mathcal{A}=(\mathcal{P},\mathcal{L},\mathcal{I}),$ has
two fixed points about an dilation then he dilation is \textit{identical
dilation} $id_{\mathcal{P}}$ of his (Theorem 2.12 in \cite{5}).

In an affine plane related to dilation $\delta \neq id_{\mathcal{P}}$ all
traces \ $P\delta \left( P\right) $ for all $P\in \mathcal{P}$ , or cross
the by a single point, or are parallel between themselves (see \cite{15}, \cite{5}).

\begin{definition} \cite{11}, \cite{5}
The \textbf{Translation} of an affine plane $\mathcal{A}=(%
\mathcal{P},\mathcal{L},\mathcal{I}),$ called the his identical dilation $id_{
\mathcal{P}}$ and every other of it's dilation, about which that the affine
plane has not fixed points.
\end{definition}

If $\sigma $ is an translation different from identical translation $id_{%
\mathcal{P}},$ then, all traces related to $\sigma $ form the
a set of parallel lines.

\begin{definition} \label{dif.direction}
For one translation $\sigma \neq id_{\mathcal{P}},$  the parallel equivalence classes of
the cleavage $\pi =\mathcal{L}/\Vert ,$ which contained tracks by  $\sigma $
of points of the plane $\mathcal{A}=(\mathcal{P},\mathcal{L},\mathcal{I})$
called the direction of his translation $\sigma $\ and marked with $\pi
_{\sigma }.$
\end{definition}

So, for $\sigma \neq id_{\mathcal{P}},$ the direction $\pi _{\sigma
}$ represented by single the trace (which is otherwise called, representative of direction) by  $\sigma $ for every point $P\in 
\mathcal{P},$ for translation $id_{\mathcal{P}}$  has
undefined direction.

Let it be $\alpha :\mathbf{Tr}_{\mathcal{A}}\longrightarrow \mathbf{Tr}_{%
\mathcal{A}},$ an whatever map of $\mathbf{Tr}_{\mathcal{A}},$ \ on
yourself. For every translation $\sigma,$ its image $\alpha \left(\sigma
\right)$ is again an translation, that can be $\alpha \left( \sigma \right)
=id_{\mathcal{P}}$ or $\alpha \left( \sigma \right) \neq id_{\mathcal{P}}$ .
So there is a certain direction or indefinite. The first equation, in
the case where \ $\sigma =id_{\mathcal{P}},$ takes the view $\alpha \left(
id_{\mathcal{P}}\right) =id_{\mathcal{P}},$ and the second $\alpha \left(
\sigma \right) \neq id_{\mathcal{P}},$ that it is not possible to $\alpha $
is map. To avoid this, yet accept that for every map $
\alpha :\mathbf{Tr}_{\mathcal{A}}\longrightarrow \mathbf{Tr}_{\mathcal{A}},$
is true this,
\[
\alpha \left( id_{\mathcal{P}}\right) =id_{\mathcal{P}}.
\]

\begin{proposition} \cite{5}
If translations $\sigma _{1}$ and $\sigma _{2}$ have the same direction with
translation $\sigma $ to an affine plane $\mathcal{A}=(\mathcal{P}
,\mathcal{L},\mathcal{I})$, then and composition $%
\sigma _{2}\circ \sigma _{1}$ has the same the direction, otherwise%
\begin{equation} \label{13}
\forall \sigma _{1},\sigma _{2},\sigma \in \mathbf{Tr}_{\mathcal{A}},\pi
_{\sigma _{1}}=\pi _{\sigma _{2}}=\pi _{\sigma }\Longrightarrow \pi _{\sigma
_{2}\circ \sigma _{1}}=\pi _{\sigma }.
\end{equation}
\end{proposition}

\begin{theorem} \cite{5}
Set $\mathbf{Tr}_{\mathcal{A}}$ of translations to an affine plane $\mathcal{%
A}$ form a group about the composition $\circ ,$ which is a sub-group of the
group $\left( \mathbf{Dil}_{\mathcal{A}},\circ \right) $ to dilations of
affine plane $\mathcal{A}$.
\end{theorem}

\begin{theorem} \cite{5} \label{thm.normal}
Group \ $\left( \mathbf{Tr}_{\mathcal{A}},\circ \right) $ of translations to
the affine plane $\mathcal{A}$ is normal sub- group of the group of
dilation’s \ $\left( \mathbf{Dil}_{\mathcal{A}},\circ \right) $ of him
plane.
\[
\forall \delta \in Dil_{\mathcal{A}},\forall \sigma \in Tr_{\mathcal{A}}\Rightarrow {\delta}^{-1} \circ \sigma \circ \delta \in Tr_{\mathcal{A}}.
\]
\end{theorem}

\begin{corollary} \cite{5} \label{cor.1}
For every dilation $\delta \in \mathbf{Dil}_{\mathcal{A}}$ and for every
translations $\sigma \in \mathbf{Tr}_{\mathcal{A}}$ of affine plane $%
\mathcal{A}\emph{=(}\mathcal{P}\emph{,}\mathcal{L}\emph{,}\mathcal{I}\emph{)}
$, translations $\sigma $ and $\delta ^{-1}\circ \sigma \circ \delta $ of
his have the same direction.
\end{corollary}

\begin{corollary} \label{14}
The translations group $\left( \mathbf{Tr}_{\mathcal{A}},\circ \right) $
of an affine plane $\mathcal{A}$ is (Abelian) commutative Group.
\end{corollary}

\section{The Endomorphisms Algebra, of the Translation's Group and their trace-preserving Associative Unitary Ring}

Consider the affine plane $\mathcal{A}=(\mathcal{P},\mathcal{L},\mathcal{I})$ and set of maps, of commutative Group $(Tr_{\mathcal{A}}, \circ) $ of affine
plane $\mathcal{A}$ in itself (see \cite{5}), so ${(Tr_{\mathcal{A}})}^{Tr_{\mathcal{A}}}=\left\{\alpha | \alpha: Tr_{\mathcal{A}}\rightarrow Tr_{\mathcal{A}}\right\}$. Let be $\alpha, \beta$ two different maps, such. Then for every $\sigma \in Tr_{\mathcal{A}}  \Rightarrow \alpha(\sigma) \in Tr_{\mathcal{A}}, \beta(\sigma) \in Tr_{\mathcal{A}}$ and $[\alpha \circ \beta](\sigma)=\alpha (\beta(\sigma))\in Tr_{\mathcal{A}}$. 
From the latter it turns out that the action of composition '$\circ$' in the set ${(Tr_{\mathcal{A}})}^{Tr_{\mathcal{A}}}$, is action induced by the action of composition '$\circ$' in the set $({\mathcal{P}})^{\mathcal{P}}$, of maps ${\mathcal{P}}\rightarrow {\mathcal{P}}$ of the affine plane $\mathcal{A}$. If the associate translations $\sigma$, the unique translation  $\alpha (\sigma) \circ \beta(\sigma)$, obtained a new map $Tr_{\mathcal{A}}\rightarrow Tr_{\mathcal{A}}$, which we call addition of $\alpha$ with $\beta$, and mark with $\alpha + \beta$. 

\begin{definition} {\label{Def.Sum}}
For every two maps  $\alpha, \beta \in {(Tr_{\mathcal{A}})}^{Tr_{\mathcal{A}}}$, the \textbf{addition} of them, that is marked $\alpha + \beta$, is called the map '${\alpha + \beta}:Tr_{\mathcal{A}} \rightarrow Tr_{\mathcal{A}}$', defined by,
\begin{equation} \label{15}
(\alpha + \beta)(\sigma)=(\alpha)(\sigma)\circ (\beta)(\sigma), \forall \sigma \in Tr_{\mathcal{A}}.
\end{equation}
\end{definition}

Accompanying any two maps $\alpha, \beta$ their sum $\alpha + \beta$, we obtain a new binary action in $Tr_{\mathcal{A}}$, that we call the addition of maps, of translations in the affine plane $\mathcal{A}$.

Thus obtained, algebra with two binary operations $({(Tr_{\mathcal{A}})}^{Tr_{\mathcal{A}}},+, \circ)$, where the sum of the two elements, whatsoever, its $\alpha, \beta$ in ${(Tr_{\mathcal{A}})}^{Tr_{\mathcal{A}}}$, is given by Definition \ref{Def.Sum}, and their composition is given by,

\begin{equation} \label{16}
[\alpha \circ \beta](\sigma)=\alpha (\beta(\sigma)), \forall \sigma \in Tr_{\mathcal{A}}
\end{equation}

\begin{definition}
The Algebra $\left({(Tr_{\mathcal{A}})}^{Tr_{\mathcal{A}}},+, \circ\right)$, is called the algebra of, maps of $Tr_{\mathcal{A}}$ on himself.
\end{definition}

A map $\alpha:Tr_{\mathcal{A}} \longrightarrow Tr_{\mathcal{A}}$, is an endomorphism of the group $(Tr_{\mathcal{A}}, \circ)$, on himself (see \cite{3}, \cite{15}), namely such that,
\begin{equation} \label{17}
\forall \sigma_{1}, \sigma_{2} \in Tr_{\mathcal{A}}, \alpha (\sigma_{1} \circ \sigma_{2})=\alpha (\sigma_{1}) \circ \alpha (\sigma_{2}).
\end{equation}

\begin{lemma} \label{lem.1}
The addition of, each two endomorphisms of $Tr_{\mathcal{A}}$ on himself, is a endomorphisms of $Tr_{\mathcal{A}}$ on himself.
\end{lemma}

\proof Let them be $\alpha, \beta$ two endomorphisms of $Tr_{\mathcal{A}}$ on himself. According to eq. \ref{17}, $\forall \sigma_{1}, \sigma_{2} \in Tr_{\mathcal{A}}$ have,
\[
\alpha (\sigma_{1} \circ \sigma_{2})=\alpha (\sigma_{1}) \circ \alpha (\sigma_{2}) \text{ \ and \  }  \beta (\sigma_{1} \circ \sigma_{2})=\beta (\sigma_{1}) \circ \beta (\sigma_{2}).
\]
Then, we have:
\begin{align*}
[\alpha + \beta](\sigma_{1} \circ \sigma_{2}) &= \alpha (\sigma_{1} \circ \sigma_{2}) \circ \beta(\sigma_{1} \circ \sigma_{2}) \ (by, \ eq. \ref{15}) \\
&= [\alpha (\sigma_{1}) \circ \alpha (\sigma_{2})] \circ [\beta (\sigma_{1}) \circ \beta (\sigma_{2})] \ (by, \ eq. \ref{17})\\
&= [\alpha (\sigma_{1}) \circ \beta (\sigma_{1})]  \circ [\alpha (\sigma_{2})]  \circ \beta (\sigma_{2})] \ (by, \ eq. \ref{15})\\
&= [\alpha + \beta](\sigma_{1}) \circ [\alpha + \beta](\sigma_{2}) \ (by, \ eq. \ref{15}).
\end{align*}
Thus,
\[
\forall \sigma_{1}, \sigma_{2} \in Tr_{\mathcal{A}},  \]
\[ 
[\alpha + \beta](\sigma_{1} \circ \sigma_{2}) = [\alpha + \beta](\sigma_{1}) \circ [\alpha + \beta](\sigma_{2}).
\]

\qed

\begin{lemma} \label{lem.2}
The composition of any two endomorphisms of $Tr_{\mathcal{A}}$ on himself, is an endomorphisms of $Tr_{\mathcal{A}}$ on himself.
\end{lemma}

\proof In the conditions when $\alpha$ and $\beta$ are endomorphisms, by eq.\ref{16} we have too 
\begin{align*}
[\alpha \circ \beta](\sigma_{1} \circ \sigma_{2}) &= \alpha [\beta (\sigma_{1} \circ \sigma_{2})] \ (by, \ eq. \ref{16}) \\
&= \alpha [\beta (\sigma_{1}) \circ \beta (\sigma_{2})] \ (by, \ eq. \ref{17})\\
&= \alpha [\beta (\sigma_{1})]  \circ \alpha [\beta (\sigma_{2})]  \  (by, eq. \ref{17})\\
&= [\alpha \circ  \beta](\sigma_{1}) \circ [\alpha \circ \beta](\sigma_{2}) \ (by, \ eq. \ref{16})
\end{align*}

Thus,
\[
\forall \sigma_{1}, \sigma_{2} \in Tr_{\mathcal{A}},
\]
\[  
[\alpha \circ \beta](\sigma_{1} \circ \sigma_{2}) = [\alpha \circ  \beta](\sigma_{1}) \circ [\alpha \circ \beta](\sigma_{2}).
\]

\qed

From Lemmas ~\ref{lem.1} \& ~\ref{lem.2}, based on the understanding of a substructure of an algebraic structure (see \cite{Lang2002}, \cite{Rotman2010}, \cite{RobertWisbauer.RingTheory1991}), we get this too,

\begin{theorem} \label{thm.1}
The set of endomorphisms of $Tr_{\mathcal{A}}$ on himself, regarding actions 'addition $+$' and 'composition $\circ$' in ti, is a substructure of algebra $\left({(Tr_{\mathcal{A}})}^{Tr_{\mathcal{A}}},+, \circ\right)$ of maps, of $Tr_{\mathcal{A}}$ on oneself.
\end{theorem}

We call it, the '\textit{endomorphisms-algebra}' of the $Tr_{\mathcal{A}}$ group on ourselves and mark with $End_{Tr_{\mathcal{A}}}$.


\begin{definition}
The '\textbf{tracer-preserving}' endomorphism of the group $(Tr_{\mathcal{A}}, \circ)$ above itself, it is called an endomorphism $\alpha \in End_{Tr_{\mathcal{A}}}$ his, such that
\begin{equation} \label{18}
\forall \sigma \in Tr_{\mathcal{A}}, \pi_{\alpha(\sigma)}=\pi_{\sigma},
\end{equation}
otherwise, any trace according to $\alpha(\sigma)$ is a trace according to $\sigma$.
\end{definition}
The 'tracer-preserving' endomorphism of the group $(Tr_{\mathcal{A}}, \circ)$ above itself, we will marked with, $End_{Tr_{\mathcal{A}}}^{TP}$

The map $0_{Tr_{\mathcal{A}}}: Tr_{\mathcal{A}} \longrightarrow Tr_{\mathcal{A}}$, determined by

\begin{equation} \label{19}
0_{Tr_{\mathcal{A}}}(\sigma)=id_{\mathcal{P}}, \forall \sigma \in Tr_{\mathcal{A}}.
\end{equation}
is an endomorphism of the translation group $Tr_{\mathcal{A}}$ on himself because,
\[
\forall \sigma_{1}, \sigma_{2} \in Tr_{\mathcal{A}}, 
\]
\[
0_{Tr_{\mathcal{A}}}(\sigma_{1} \circ \sigma_{2})=id_{\mathcal{P}}=id_{\mathcal{P}} \circ id_{\mathcal{P}} = 0_{Tr_{\mathcal{A}}}(\sigma_{1}) \circ 0_{Tr_{\mathcal{A}}}(\sigma_{2}).
\]
So,
\[
0_{Tr_{\mathcal{A}}} \in End_{Tr_{\mathcal{A}}}
\]
The endomorphism $0_{Tr_{\mathcal{A}}}$ call it, \textbf{zero endomorphism} of $Tr_{\mathcal{A}}$ on himself.

From eq.\ref{19} also indicate that, 
\[
\forall \sigma \in Tr_{\mathcal{A}},  \pi_{0_{Tr_{\mathcal{A}}}(\sigma)}=\pi_{\sigma}.
\]
so we have that,
\[
0_{Tr_{\mathcal{A}}} \in End_{Tr_{\mathcal{A}}}^{TP},
\]
thus is i true this,

\begin{proposition}
The zero endomorphism $0_{Tr_{\mathcal{A}}}$ of $Tr_{\mathcal{A}}$ on himself, is a trace-preserving endomorphism of $Tr_{\mathcal{A}}$ on himself. 
\end{proposition}

The identical map $1_{Tr_{\mathcal{A}}}: Tr_{\mathcal{A}} \longrightarrow Tr_{\mathcal{A}}$, defined by
\begin{equation} \label{20}
1_{Tr_{\mathcal{A}}}(\sigma)=\sigma, \forall \sigma \in Tr_{\mathcal{A}},
\end{equation}
is also, an endomorphism of the translation group $Tr_{\mathcal{A}}$ on himself because,
\[
\forall \sigma_{1}, \sigma_{2} \in Tr_{\mathcal{A}}, 
\]
\[
1_{Tr_{\mathcal{A}}}(\sigma_{1} \circ \sigma_{2})=\sigma_{1} \circ \sigma_{2}= 1_{Tr_{\mathcal{A}}}(\sigma_{1}) \circ 1_{Tr_{\mathcal{A}}}(\sigma_{2}).
\]
So,
\[
1_{Tr_{\mathcal{A}}} \in End_{Tr_{\mathcal{A}}}
\]
The endomorphism $1_{Tr_{\mathcal{A}}}$ call it, \textbf{unitary endomorphism} of $Tr_{\mathcal{A}}$ on himself.
From \ref{20} we have immediately that
\[
\forall \sigma \in Tr_{\mathcal{A}},  \pi_{1_{Tr_{\mathcal{A}}}(\sigma)}=\pi_{\sigma},
\]
so we have that,
\[
1_{Tr_{\mathcal{A}}} \in End_{Tr_{\mathcal{A}}}^{TP},
\]
thus is i true this,
\begin{proposition}
The unitary endomorphism '$1_{Tr_{\mathcal{A}}}$' of $Tr_{\mathcal{A}}$ on himself, is a trace-preserving endomorphism of $Tr_{\mathcal{A}}$ on himself.
\end{proposition}

\begin{theorem} \label{thm.2}
If $\alpha$ and $\beta$ are two trace-preserving endomorphisms of $Tr_{\mathcal{A}}$ on himself, then their sum $\alpha + \beta$ is an trace-preserving endomorphism.
\end{theorem}

\proof In terms of Theorem, according to Lemma ~\ref{lem.1}, the addition $\alpha + \beta$ is an endomorphisms of $Tr_{\mathcal{A}}$ on himself. We note $\sigma_{1}=\alpha(\sigma)$ and $\sigma_{2}=\beta(\sigma)$. Given that $\alpha$ and $\beta$ are two trace-preserving endomorphisms, according to eq. \ref{18}, $\forall \sigma \in Tr_{\mathcal{A}}$ have,
\begin{equation} \label{18'}
\pi_{\alpha(\sigma)}=\pi_{\sigma_{1}}=\pi_{\sigma} \text{ \ and \ } \pi_{\beta(\sigma)}=\pi_{\sigma_{2}}=\pi_{\sigma}
\end{equation}
From here, according to \ref{13} we get
\begin{equation} \label{13'}
\forall \sigma \in Tr_{\mathcal{A}}, \pi_{\alpha(\sigma) \circ \beta(\sigma)}=\pi_{{\sigma_{1}} \circ {\sigma_{2}}}=\pi_{\sigma}
\end{equation}
Then, by eq. \ref{15} and eq. \ref{13'}, have that
\[
\pi_{(\alpha + \beta)(\sigma)} =\pi_{\alpha(\sigma) \circ \beta(\sigma)}=\pi_{\sigma}
\]
Hence,
\[
\forall \sigma \in Tr_{\mathcal{A}}, \pi_{(\alpha + \beta)(\sigma)} =\pi_{\sigma}
\]
\qed

\begin{theorem} \label{thm.3}
If $\alpha$ and $\beta$ are two trace-preserving endomorphisms of $Tr_{\mathcal{A}}$ on himself, then their composition $\alpha \circ \beta$ is an trace-preserving endomorphism.
\end{theorem}

\proof In terms of Theorem, according to Lemma ~\ref{lem.2}, the composition $\alpha \circ \beta$ is an endomorphisms of $Tr_{\mathcal{A}}$ on himself. We note $\sigma_{1}=\beta(\sigma)$, given that $\alpha$ and $\beta$ are two trace-preserving endomorphisms, according to eq. \ref{18}, $\forall \sigma \in Tr_{\mathcal{A}}$ and for $\sigma_{1}=\beta(\sigma)$, have
\begin{equation} \label{eq.18''}
\pi_{\alpha(\sigma_{1})}=\pi_{\sigma_{1}} \text{ \ and \ } \pi_{\beta(\sigma)}=\pi_{\sigma}
\end{equation}
Then, by eq. \ref{16} and eq. \ref{eq.18''}, have that
\[
\pi_{(\alpha \circ \beta)(\sigma)}=\pi_{\alpha(\beta(\sigma))}=\pi_{\alpha(\sigma_{1})}=\pi_{\sigma_{1}}=\pi_{\beta(\sigma)}=\pi_{\sigma}.
\]
Hence,
\[
\forall \sigma \in Tr_{\mathcal{A}}, \pi_{(\alpha \circ \beta)(\sigma)}=\pi_{\sigma}.
\]
\qed

We note now, with $\varphi:Tr_{\mathcal{A}}\longrightarrow Tr_{\mathcal{A}}$, defined by
\begin{equation} \label{21}
\forall \sigma \in Tr_{\mathcal{A}}, \varphi(\sigma)={\sigma}^{-1}
\end{equation}
is an endomorphism of the commutative group of translations $Tr_{\mathcal{A}}$ on himself, because $\forall \sigma_{1}, \sigma_{2} \in Tr_{\mathcal{A}},$ we have
\begin{align*}
\varphi(\sigma_{1} \circ \sigma_{2}) &= {(\sigma_{1} \circ \sigma_{2})}^{-1} \ (by, \ eq. \ref{21}) \\
&= {(\sigma_{2} )}^{-1} \circ {(\sigma_{1})}^{-1}  \\
&= {(\sigma_{1} )}^{-1} \circ {(\sigma_{2})}^{-1}  \  (by, \ref{14})\\
&= \varphi(\sigma_{1}) \circ \varphi(\sigma_{2}) 
\end{align*}
Also, we have that,
\[
\forall \sigma \in Tr_{\mathcal{A}}, \pi_{\sigma}=\pi_{{\sigma}^{-1}} \Rightarrow \pi_{\varphi(\sigma)}=\pi_{{\sigma}^{-1}}= \pi_{\sigma}.
\]
from this, we have prove this
\begin{proposition}
The endomorphisms $\varphi:Tr_{\mathcal{A}}\longrightarrow Tr_{\mathcal{A}}$, defined by $\varphi(\sigma)={\sigma}^{-1},\forall \sigma \in Tr_{\mathcal{A}},$ is an trace-preserving endomorphisms of $Tr_{\mathcal{A}}$ on himself.
\end{proposition}

\begin{proposition}
For a trace-preserving endomorphism $\alpha \in End_{Tr_{\mathcal{A}}}^{TP}$, the map $-\alpha: Tr_{\mathcal{A}}\longrightarrow Tr_{\mathcal{A}}$, defined by,
\begin{equation} \label{21'}
\forall \sigma \in Tr_{\mathcal{A}}, (-\alpha)(\sigma)=(\alpha(\sigma))^{-1} \Rightarrow (-\alpha) \in End_{Tr_{\mathcal{A}}}^{TP}
\end{equation}
well $-\alpha$ is an trace-preserving endomorphism of $Tr_{\mathcal{A}}$ on himself.
\end{proposition}

\proof Note that $\forall \sigma \in Tr_{\mathcal{A}}$,
\begin{align*}
(-\alpha)(\sigma) &= (\alpha(\sigma))^{-1}  \\
&= \varphi(\alpha(\sigma)) \  (by, eq. \ref{21}) \\
&= (\varphi \circ \alpha)(\sigma) \  (by, eq. \ref{16}) 
\end{align*}
which indicates that, $-\alpha=\varphi \circ \alpha$ and according to Theorem \ref{thm.3}, he is a  trace-preserving endomorphism.
\qed

The endomorphism $'-\alpha'$ we call it the additive inverse endomorphism of endomorphism $\alpha$.
Consider now the set 
\begin{equation} \label{23}
End_{Tr_{\mathcal{A}}}^{TP}=\left\{\alpha \in End_{Tr_{\mathcal{A}}} \ | \ \alpha - \text{is a trace-preserving endomrphism}\right\}
\end{equation}
of trace-preserving endomorphisms of $Tr_{\mathcal{A}}$ in itself. According to Theorem \ref{thm.2} and Theorem \ref{thm.3}, $(End_{\mathcal{A}}^{TP}, +, \circ)$ is a substructure of the algebra  $(End_{\mathcal{A}}, +, \circ)$, of endomorphisms related to addition and composition actions, therefore it is itself an algebra.

\begin{theorem} \label{additive.group}
The Grupoid $(End_{Tr_{\mathcal{A}}}^{TP}, +)$, is commutative (Abelian) Group.
\end{theorem}

\proof 1. $\forall \alpha, \beta, \gamma \in End_{Tr_{\mathcal{A}}}^{TP}$, have
\begin{align*}
\forall \sigma \in Tr_{\mathcal{A}}, [(\alpha + \beta)+\gamma](\sigma) &= (\alpha + \beta)(\sigma) \circ (\gamma)(\sigma) \ (by, \ eq. \ref{15}) \\
&= [\alpha(\sigma) \circ (\beta)(\sigma)] \circ (\gamma)(\sigma) \ (by, \ eq. \ref{15}) \\
&=\alpha(\sigma) \circ [(\beta)(\sigma) \circ (\gamma)(\sigma) ] \ (by, \   \ref{14})\\
&=\alpha(\sigma) \circ [(\beta + \gamma)(\sigma) ] \ (by, \ eq. \ref{15})\\
&= [\alpha + (\beta+\gamma)](\sigma)  \ (by, \ eq. \ref{15})
\end{align*}
So,
\[
\forall \alpha, \beta, \gamma \in End_{Tr_{\mathcal{A}}}^{TP}, (\alpha + \beta)+\gamma= \alpha + (\beta+\gamma).
\]

2. $\forall \alpha \in End_{Tr_{\mathcal{A}}}^{TP}$, have
\begin{align*}
\forall \sigma \in Tr_{\mathcal{A}}, [\alpha + 0_{Tr_{\mathcal{A}}}](\sigma) &= \alpha (\sigma) \circ 0_{Tr_{\mathcal{A}}}(\sigma) \ (by, \ eq. \ref{15}) \\
&= \alpha (\sigma) \circ id_{\mathcal{P}} \ (by, \ eq. \ref{19}) \\
&= \alpha (\sigma)
\end{align*}
So exists, the zero element in $End_{Tr_{\mathcal{A}}}^{TP}$, which is zero endomorphism $0_{Tr_{\mathcal{A}}}$ for which, we have proven that it is a trace-preserving endomorphism, and we have that, 
\[
\forall \alpha \in End_{Tr_{\mathcal{A}}}^{TP}, \exists 0_{Tr_{\mathcal{A}}}\in End_{Tr_{\mathcal{A}}}^{TP}, \alpha + 0_{Tr_{\mathcal{A}}}= 0_{Tr_{\mathcal{A}}} +\alpha =\alpha.
\]

3. $\forall \alpha \in End_{Tr_{\mathcal{A}}}^{TP}$, have
\begin{align*}
\forall \sigma \in Tr_{\mathcal{A}}, [\alpha + (-\alpha)](\sigma) &= \alpha (\sigma) \circ (-\alpha)(\sigma) \ (by, \ eq. \ref{15}) \\
&= \alpha (\sigma) \circ [\alpha (\sigma)]^{-1} \ (by, \ eq. \ref{21'}) \\
&= id_{\mathcal{P}} \\
&= 0_{Tr_{\mathcal{A}}}. \ (by, \ eq. \ref{19})
\end{align*}
So exists, the additive inverse element of $\alpha$ in $End_{Tr_{\mathcal{A}}}^{TP}$, which is $-\alpha$ for which, we have proven that it is a trace-preserving endomorphism, and we have that, 
\[
\forall \alpha \in End_{Tr_{\mathcal{A}}}^{TP}, \exists (-\alpha) \in End_{Tr_{\mathcal{A}}}^{TP}, \alpha + (-\alpha)= 0_{Tr_{\mathcal{A}}}.
\]

4. $\forall \alpha , \beta \in End_{Tr_{\mathcal{A}}}^{TP}$, have
\begin{align*}
\forall \sigma \in Tr_{\mathcal{A}}, [\alpha + \beta](\sigma) &= \alpha (\sigma) \circ \beta (\sigma) \ (by, \ eq. \ref{15}) \\
&= \beta (\sigma) \circ \alpha (\sigma) \ (by, \  \ref{14}) \\
&= [\beta + \alpha] (\sigma). \ (by, \ eq. \ref{15})
\end{align*}
So,
\[
\forall \alpha , \beta \in End_{Tr_{\mathcal{A}}}^{TP}, \alpha + \beta=\beta + \alpha.
\]
\qed

\begin{theorem} \label{unitary.ring}
The algebra $(End_{Tr_{\mathcal{A}}}^{TP}, + , \circ)$, is a associative unitary Ring.
\end{theorem}

\proof According to the Definition of a associative Unitary Ring \cite{Lang2002}, \cite{Rotman2010}, \cite{RobertWisbauer.RingTheory1991} we are required to prove the following conditions, \\
(1) The grupoid $(End_{Tr_{\mathcal{A}}}^{TP}, + )$, is commutative group, we have proved this in the Theorem \ref{additive.group}. \\
(2) The action '$\circ$' is associative in $End_{Tr_{\mathcal{A}}}^{TP}$, truly, by definition of composition, we have
\begin{align*}
\forall \sigma \in Tr_{\mathcal{A}}, [(\alpha \circ \beta) \circ \gamma] (\sigma) &= (\alpha \circ \beta)(\gamma (\sigma))  \\
&= \alpha ( \beta (\gamma (\sigma)))  \\
&= \alpha [(\beta \circ \gamma) (\sigma)] \\
&= [\alpha \circ( \beta \circ \gamma)] (\sigma).
\end{align*}
So we have that,
\[
\forall \alpha, \beta, \gamma \in End_{Tr_{\mathcal{A}}}^{TP}, (\alpha \circ \beta) \circ \gamma=\alpha \circ( \beta \circ \gamma).
\]

(3) The composition is 'distributive' according to 'addition', ie,
\[
\forall \alpha, \beta, \gamma \in End_{Tr_{\mathcal{A}}}^{TP}, \]
\[
\alpha \circ ( \beta + \gamma ) =\alpha \circ  \beta + \alpha \circ \gamma , \\
\text{ \ and \ } (\alpha + \beta) \circ \gamma=\alpha \circ \gamma + \beta \circ \gamma.
\]
Really,
\begin{align*}
\forall \sigma \in Tr_{\mathcal{A}}, [\alpha \circ ( \beta + \gamma )] (\sigma) &= \alpha  [( \beta + \gamma )(\sigma)] \ \ (by, \ eq. \ref{16})  \\
&= \alpha [\beta(\sigma) \circ \gamma (\sigma)] \ \ (by, \ eq. \ref{15})  \\
&= \alpha [\beta(\sigma)] \circ \alpha [\gamma (\sigma)] \ \ (by, \ eq. \ref{17})  \\
&= [\alpha \circ \beta](\sigma) \circ [\alpha \circ \gamma] (\sigma) \ \ (by, \ eq. \ref{16})  \\
&= [\alpha \circ \beta + \alpha \circ \gamma] (\sigma). \ \ (by, \ eq. \ref{15})  
\end{align*}

Hence
\[ \forall \alpha, \beta, \gamma \in End_{Tr_{\mathcal{A}}}^{TP}, \alpha \circ ( \beta + \gamma ) =\alpha \circ  \beta + \alpha \circ \gamma. \]

Also,
\begin{align*}
\forall \sigma \in Tr_{\mathcal{A}}, [(\alpha + \beta) \circ \gamma] (\sigma) &= [\alpha + \beta] (\gamma (\sigma)) \ \ (by, \ eq. \ref{16})  \\
&= \alpha (\gamma (\sigma))  \circ  \beta (\gamma (\sigma))  \ \ (by, \ eq. \ref{15})  \\
&= [\alpha \circ \gamma] (\sigma)  \circ  [\beta \circ \gamma] (\sigma)  \ \ (by, \ eq. \ref{16})  \\
&= [\alpha \circ \gamma+ \beta \circ \gamma] (\sigma). \ \ (by, \ eq. \ref{15})  
\end{align*}

Hence
\[ \forall \alpha, \beta, \gamma \in End_{Tr_{\mathcal{A}}}^{TP}, (\alpha + \beta) \circ \gamma=\alpha \circ \gamma + \beta \circ \gamma.
\]

(4) In $End_{Tr_{\mathcal{A}}}^{TP}$, exist the unitary element, related to composition
\[
\forall \alpha \in End_{Tr_{\mathcal{A}}}^{TP}, \text{ \ have,}
\]
\begin{align*}
\forall \sigma \in Tr_{\mathcal{A}}, [\alpha \circ 1_{Tr_{\mathcal{A}}}] (\sigma) &=\alpha [ 1_{Tr_{\mathcal{A}}} (\sigma)] \ \ (by, \ eq. \ref{16})  \\
&= \alpha (\sigma)   \ \ (by, \ eq. \ref{20}). 
\end{align*}

Hence

\[ \forall \alpha \in End_{Tr_{\mathcal{A}}}^{TP}, \alpha \circ 1_{Tr_{\mathcal{A}}}=\alpha,
\]

ie, the unitary element of $End_{Tr_{\mathcal{A}}}^{TP}$, is the unitary endomorphism '$1_{Tr_{\mathcal{A}}}$' of $Tr_{\mathcal{A}}$ on himself.

\qed


\bibliographystyle{amsplain}
\bibliography{NSrefs}

\end{document}